\documentclass[11pt]{amsart}
\usepackage{enumerate}
\usepackage{geometry}                
\usepackage{graphicx}
\usepackage{amssymb}
\usepackage{epstopdf}
\usepackage{geometry} 
\usepackage{graphicx}
\usepackage{xfrac}
\usepackage{amsmath} 
\usepackage{stackengine}
\usepackage{romanbar}
\numberwithin{equation}{section}
\usepackage{commath}
\usepackage{cleveref}
\usepackage[normalem]{ulem}
\usepackage{bbm}
\usepackage{subfig}
\usepackage{float}
\usepackage{url}
\usepackage{mathtools}
\usepackage{amsthm} 
\usepackage[disable]{todonotes} 
\newcommand{\R}{\mathbb R}

\theoremstyle{plain}
\newtheorem{theorem}{Theorem}[section]
\newtheorem*{theorem*}{Theorem}
\newtheorem{lemma}[theorem]{Lemma}

\newtheorem{proposition}[theorem]{Proposition}
\newtheorem{conjecture}[theorem]{Conjecture}

\theoremstyle{remark}
\newtheorem{remark}{Remark}
\theoremstyle{remark}

\theoremstyle{definition}
\newtheorem{definition}[theorem]{Definition}
\allowdisplaybreaks
\title[(Non)uniqueness of contact discontinuities for compressible Euler]{Contact discontinuities for 2-D isentropic Euler are unique in 1-D but wildly non-unique otherwise}
\author[Sam G. Krupa]{Sam G.  Krupa}
\address[Sam G. Krupa]{\newline Max Planck Institute for Mathematics in the Sciences, \newline 04103 Leipzig, Germany}
\email{Sam.Krupa@mis.mpg.de}

\author[L\'{a}szl\'{o} Sz\'{e}kelyhidi, Jr.]{L\'{a}szl\'{o} Sz\'{e}kelyhidi, Jr.}
\address[L\'{a}szl\'{o} Sz\'{e}kelyhidi, Jr.]{\newline Max Planck Institute for Mathematics in the Sciences, \newline 04103 Leipzig, Germany}
\email{Laszlo.Szekelyhidi@mis.mpg.de}

\thanks{Funded by the German Research Foundation (DFG) project number 525859002. Part of this work was completed while the first author was visiting National Tsing Hua University (Taiwan) and the National Center for Theoretical Sciences (Taiwan). The first author would like to thank his host Jin-Cheng Jiang at National Tsing Hua University and the National Center for Theoretical Sciences for the wonderful working atmosphere.}

\date{\today}                                           

\begin{document}
\keywords{Systems of conservation laws, contact discontinuity, vortex sheet, compressible Euler, isentropic Euler, non-uniqueness, convex integration, multiple space dimensions.}
\subjclass[2020]{Primary 35L65; Secondary 35A02, 35B30, 35L45, 35Q31, 76N15}
\begin{abstract}
   We develop a general framework for studying non-uniqueness of the Riemann problem for the isentropic compressible Euler system in two spatial dimensions, and in this paper we present the most delicate result of our method: non-uniqueness of the contact discontinuity. Our approach is computational, and uses the pressure law as an additional degree of freedom. 
   
    The stability of the contact discontinuities for this system is a major open problem (see Gui-Qiang Chen and Ya-Guang Wang [{\em Nonlinear partial differential equations}, volume~7 of {\em Abel
  Symposia}. Springer, Heidelberg, 2012.]). 
  
  We find a smooth pressure law $p$, verifying the physically relevant condition $p'>0$, such that for the isentropic compressible Euler system with this pressure law, contact discontinuity initial data is wildly non-unique in the class of bounded, admissible weak solutions. This result resolves the question of uniqueness for contact discontinuity solutions in the compressible regime. 
  
  Moreover, in the \emph{same regularity class} in which we have non-uniqueness of the contact discontinuity, i.e. $L^\infty$, with no $BV$ regularity or self-similarity, we show that the classical contact discontinuity solution to the two-dimensional isentropic compressible Euler system is in fact \emph{unique} in the class of bounded, admissible weak solutions if we restrict to 1-D solutions.
\end{abstract}

\maketitle

\section{Introduction}
This paper focuses on the initial value problem for the isentropic compressible Euler system 
\begin{equation}
\begin{cases}\label{system}
\partial_t \rho + \text{div}_x (\rho v) = 0  \\
\partial_t (\rho v) + \text{div}_x (\rho v \otimes v) + \nabla_x [p(\rho)] = 0  \\
\rho(\cdot, 0) = \rho^0  \\
v(\cdot, 0) = v^0. 
\end{cases}
\end{equation}
We write $\rho$ for the density and $v$ for the velocity of the fluid. The pressure $p$ is some given function of $\rho$. We say the system \eqref{system} is \emph{hyperbolic} when $p' > 0$. This is the case which is relevant for physical models arising from continuum mechanics (see e.g. \cite[p.~56]{MR3289359}). We consider the two-dimensional setting in full space, i.e.~with spatial domain $\R^2$, and denote the space variables as $x = (x_1, x_2) \in \mathbb{R}^2$. Similarly, we write $v = (v_1, v_2) \in \mathbb{R}^2$ for the components of the vector field $v=v(x,t)$.

In this paper we will consider \emph{admissible weak solutions}. These are pairs $\rho,v\in L^\infty(\R^2\times (0,T)$ which verify \eqref{system} in the sense of distributions, and in addition satisfy the \emph{entropy inequality}. In the strong form the latter reads
\begin{align}
\partial_t \left( \rho \epsilon(\rho) + \rho \frac{|v|^2}{2} \right) + \text{div}_x \left[ \left( \rho \epsilon(\rho) + \rho \frac{|v|^2}{2} + p(\rho) \right) v \right] \leq 0, \label{system_adm}
\end{align}
where the internal energy density $\epsilon(\rho)$ is related to the pressure through the relation 
\begin{align}\label{pressure_relation}
    p(\rho) = \rho^2 \epsilon'(\rho).
\end{align}
For weak solutions \eqref{system_adm} is also to be understood in the sense of distributions, i.e.
\begin{multline}\label{system_admin_integral}
    \int_0^\infty \int_{\mathbb{R}^2} \left[ \left( \rho \epsilon(\rho) + \rho \frac{|v|^2}{2} \right) \partial_t \varphi + \left( \rho \epsilon(\rho) + \rho \frac{|v|^2}{2} + p(\rho) \right) v \cdot \nabla_x \varphi \right] dx dt
\\
+ \int_{\mathbb{R}^2} \left( \rho^0 (x) \epsilon(\rho^0 (x)) + \rho^0 (x) \frac{|v^0 (x)|^2}{2} \right) \varphi (x, 0) dx \geq 0,
\end{multline}
for every nonnegative test function $\varphi \in C_c^\infty (\mathbb{R}^2 \times [0, \infty))$.

In the case of the two-dimensional isentropic Euler equations, the \emph{Riemann problem} concerns a specific class of initial data of the form
\begin{align}
(\rho^0 (x), v^0 (x)) := \begin{cases}
(\rho_-, v_-) & \text{if } x_2 < 0 \\
(\rho_+, v_+) & \text{if } x_2 > 0,
\end{cases} \label{init_data}
\end{align}
where $\rho_{\pm}, v_{\pm} = (v_{\pm 1}, v_{\pm 2})$ are fixed constants. In the classical theory of conservation laws in one spatial dimension, the Riemann problem admits self-similar solutions consisting of a finite number of constant states separated by shocks, rarefaction waves and contact discontinuities. The proof of existence of such self-similar solutions called \emph{fan solutions}, which reduces to the analysis of an associated algebraic system, goes back to the seminal work of Riemann \cite{Riemann1860}, see also Dafermos \cite[Theorem 9.5.1]{dafermos_big_book} and Smoller \cite[Chapter 17 and Chapter 18 \S B]{MR1301779}.

One-dimensional fan solutions can be trivially extended to higher space dimensions. However, it turns out that uniqueness is lost in the two-dimensional (and higher-D) case in a quite dramatic way. This was first shown in work of Chiodaroli-De Lellis-Kreml \cite{multi_d_illposed}, building on previous work on convex integration for the incompressible and compressible Euler systems \cite{MR2564474,MR2842999}. More precisely, in \cite{multi_d_illposed} the authors showed that (i) for the pressure law $p(\rho)=\rho^2$ there exist Riemann data arising from a compression wave, for which highly non-unique class of (infinitely many!) admissible weak solutions can be constructed via convex integration, and (ii) obtained an open set of Riemann data from which such non-uniqueness arises, provided one has the freedom to choose the pressure function. It should be mentioned that the weak solutions thus constructed are admissible (in the sense of \eqref{system_adm}), but in general do not have any additional regularity properties besides being bounded and measurable, and are genuinely two-dimensional. Sometimes such weak solutions are called ``wild solutions" in the literature. 

Subsequently there has been intensive work in classifying those Riemann data for which such non-uniqueness can arise, see \cite{MR3744380,MR3269641,MR3831840,MR3816641,MR4385531,MR4210751}. For a rather comprehensive summary of the state of the art, we refer to Section 7 of the monograph \cite{MR4385531}. The conclusion of these works is that, at least with the polytropic pressure law $p(\rho)=\rho^{\gamma}$ with $\gamma>1$, whenever the classical solution to the 1-D Riemann problem contains a shock, the above non-uniqueness phenomenon persists for admissible weak solutions. 

Invariably, the technique in all these works on wild non-uniqueness involves finding a suitable subsolution to \eqref{system}-\eqref{system_adm} with given initial datum. We recall that, in general, a subsolution to \eqref{system}-\eqref{system_adm} is a triple $(\rho, v,R)\colon \mathbb{R}^2\times[0,\infty)\to (0,\infty)\times \mathbb{R}^2\times \mathcal{S}_{+}^{2 \times 2}$, $\mathcal{S}_{+}^{2 \times 2}$ denoting the set of symmetric $2 \times 2$ positive semi-definite matrices, solving the following relaxed system of equations in the sense of distributions, cf.~\cite{MR2564474}: 
 \begin{equation}
\begin{cases}\label{subsolution_system}
\partial_t \rho + \text{div}_x (\rho v) = 0,  \\
\partial_t (\rho v) + \text{div}_x (\rho v \otimes v+\rho R +p(\rho) \mathrm{Id})= 0,  \\
\partial_t \left( \rho \epsilon(\rho) + \rho \frac{|v|^2}{2} +\frac{\rho}{2}\mathrm{tr}(R) \right) + \text{div}_x \left[ \left( \rho \epsilon(\rho) + \rho \frac{|v|^2}{2} + p(\rho) +\frac{\rho}{2}\mathrm{tr}(R)\right) v \right] \leq 0.
\end{cases}
\end{equation}

Note that the above system is under-determined, because there is no equation for the Reynolds stress term $R$. In particular, a subsolution with $R\equiv 0$ is a solution of the original system \eqref{system}-\eqref{system_adm}. Conversely, if $(\rho,v,R)$ is a \emph{strict} subsolution in the sense that on some open space-time domain $R>0$ (i.e.~positive definite), then convex integration leads to wild non-uniqueness (cf.~\cite[Section 3.3]{multi_d_illposed}). Thus, in a nutshell, the existence of wild non-uniqueness for given initial datum $(v^0,\rho^0)$ boils down to the existence of a strict subsolution $(\rho,v,R)$ with
\begin{equation*}
\rho(\cdot,0)=\rho^0,\quad v(\cdot,0)=v^0,\quad R(\cdot,0)=0.	
\end{equation*}
For details we refer to \cite{MR2564474,multi_d_illposed}. For the case of the Riemann problem, i.e.~initial data of the form \eqref{init_data}, the authors in \cite{multi_d_illposed} introduced the notion of \emph{fan subsolutions}, as self-similar analogues of fan solutions depending on one space dimension, leading to a relaxation of the classical Riemann problem - see Section \ref{sec:convex_int} for details. Then, the question of existence for a fan subsolution boils down to solving a large algebraic system of equations and inequalities. 

\begin{figure}[tb]
\centering
      \includegraphics[width=\textwidth]{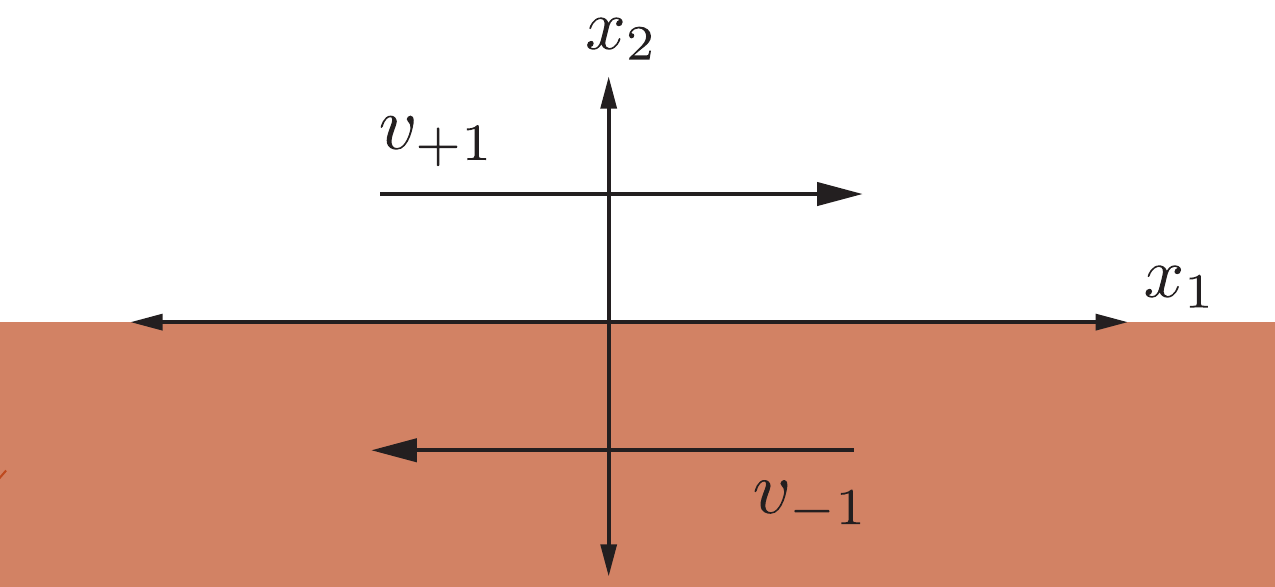}
  \caption{A diagram of the contact discontinuity initial data we consider in this paper.}\label{fig:vortex_sheet_init}
\end{figure}
 
The very low regularity of such ``wild" solutions makes their physical interpretation highly questionable. Nevertheless, as has been argued in the context of hydrodynamic turbulence \cite{MR4395517}, the significance of such non-uniqueness lies in the corresponding strong instability of the one-dimensional Riemann solution with respect to genuine multi-dimensional perturbations. Indeed, the case most studied in the incompressible setting, namely vortex sheets, is intimately related to the Kelvin-Helmholtz instability, and it is well-known, both from experiments and numerical simulations that small perturbations lead to dramatic deviations from the vortex-evolution picture by build-up of a linearly growing turbulent zone (cf. \cite{MR2842999,MR4544797,butterfly}). The corresponding compressible picture has also intensively been investigated theoretically and numerically, we refer e.g.~to \cite{MR0085036,MR0097930,MR2095445,MR2423311,MR1021640,MR0914450,MR1015065,MR3289359,MR3648106}. Nevertheless, even for the isentropic case \eqref{system}-\eqref{system_adm}, Riemann data which correspond to such Kelvin-Helmholtz instabilities and more generally, when the Riemann solution contains no shock wave, are outside the scope of the papers cited above. A particularly relevant case is Riemann datum corresponding to a single \emph{contact discontinuity}: such Riemann data has a jump  \emph{only in the velocity direction tangential to the discontinuity}, whilst both the density $\rho$ and the normal velocity component are constant across the interface. This is also known as the \emph{compressible vortex sheet} initial data. See \Cref{fig:vortex_sheet_init}.


In this paper, we develop a computer-assisted framework for finding fan subsolutions and apply it to the contact discontinuity $v_{\pm}=(\pm 1,0)$, $\rho_{\pm}=\rho$ (const), by exploiting the additional degree of freedom afforded by the pressure law. More precisely, our first main result is:

\begin{theorem}[Non-uniqueness for contact discontinuity initial data]\label{main_theorem_nonunique}
    There exists a smooth pressure function $p\colon (0,\infty)\to\mathbb{R}$ such that $p,p'>0$, and the system \eqref{system}-\eqref{system_adm} with this pressure law admits initial data in the form \eqref{init_data} with $\rho_+=\rho_-$, $v_{+ 2}=v_{- 2}$ and $v_{+ 1}\neq v_{- 1}$ which gives rise to infinitely many bounded, non-unique solutions with density uniformly bounded away from 0. 
\end{theorem}

We wish to point out that the existence of a weak solution to \eqref{system}-\eqref{system_adm} of contact discontinuity type is independent of the choice of pressure law, unlike e.g. the case of rarefaction waves. 

Complementing Theorem \ref{main_theorem_nonunique} we show that in the \emph{same regularity class} in which we have non-uniqueness of solutions with contact discontinuity initial data (i.e., bounded and measurable functions verifying \eqref{system}-\eqref{system_adm}), we have \emph{uniqueness} when we restrict to 1-D solutions which are functions only of $x_2$.

\begin{theorem}[Uniqueness of 1-D solutions]\label{main_theorem_unique}
    Consider the isentropic compressible Euler system \eqref{system} with any pressure law $p\colon (0,\infty)\to\mathbb{R}$ verifying $p,p'>0$. Consider a solution $(\rho,v_1,v_2)\in L^\infty(\mathbb{R}^2\times[0,\infty))$ to \eqref{system}-\eqref{system_adm}, where $(\rho,v_1,v_2)$ are functions only of $x_2$, $\inf \rho>0$, and $(\rho,v_1,v_2)$ has initial data in the form \eqref{init_data} verifying $\rho_+=\rho_->0$, $v_{+2}=v_{-2}$, and $v_{+ 1}\neq v_{- 1}$. Then, $(\rho,v_1,v_2)$ is in fact the classical, self-similar solution with a single contact discontinuity.
\end{theorem}

This result is a substantial generalization of \cite[Proposition 8.1]{multi_d_illposed}, where uniqueness is shown for 1-D solutions to \eqref{system}-\eqref{system_adm} with the additional hypotheses of $BV$ regularity and self-similarity. Our result does not require these rather strong regularity hypotheses, thus highlighting the fact that it is really the spatial two (or higher) dimensionality which is responsible for the wild non-uniqueness.  In particular, we note that vanishing viscosity solutions to \eqref{system}-\eqref{system_adm} with initial datum given by \eqref{init_data} are one-dimensional, but we are not aware of a priori estimates that would guarantee BV regularity in the limit in general (see \cite{MR2150387} for the small-BV case in 1-D systems). However, our result shows that indeed, \emph{the vanishing viscosity limit for initial data of contact discontinuity type is unique.} A related recent result is \cite{MR4268832}, where the authors show the stability and uniqueness of a planar contact discontinuity without shear (i.e. entropy waves, where $v_+=v_-$), for the three-dimensional full Navier–Stokes–Fourier system, in the class of vanishing dissipation limits. See also \cite{wangwei} for a very recent result on the nonlinear stability of entropy waves.

The condition in Theorem \ref{main_theorem_unique} of depending only one one spatial dimension may seem restrictive. However, it is worth pointing out that in general the question of uniqueness of solutions to conservation laws in one spatial dimension and with an associated entropy inequality, in the class $L^\infty$, is still a major open problem in the field. In particular, see Bressan's Open Problem \# 6 \cite{2023arXiv231016707B}. The nonperturbative $L^2$ stability theory (see e.g. \cite{VASSEUR2008323,MR4487515}) is a generalization of weak/strong stability which can consider discontinuous solutions with arbitrary $L^\infty$ norm, but a trace condition must still be assumed a priori. The closely related topic of convex integration for 1-D conservation laws has very recently become an area of intense study: see for example the recent results \cite{2024arXiv240702927C,MR4144350,https://doi.org/10.48550/arxiv.2208.10979,2022arXiv221114239K,2024arXiv240307784K}, and also the earlier work \cite{MR2008346}. 

\subsection{Some comments on our proofs}\label{sec:comments}

\subsubsection{\Cref{main_theorem_nonunique}: Non-uniqueness}
Our proof of \Cref{main_theorem_nonunique} uses the framework and setup of convex integration and the search for admissible fan subsolutions, as has been developed in \cite{multi_d_illposed} and used subsequently. We refer to the survey of the state of the art in \cite{MR4385531} and the exposition in Section \ref{sec:convex_int}, in particular \Cref{fan_def}. This approach starts by fixing the number of waves in the fan subsolution and then solving the corresponding algebraic system of equations and (strict) inequalities (see \Cref{RH_prop}). In the latter, it is advantageous to keep the pressure law $p$ undetermined, as this allows one to avoid the non-linearity in the algebraic system arising from the pressure (cf.~\cite[Section 7]{multi_d_illposed} and also \cite{MR2048569,2024arXiv240307784K} in related contexts).    In most examples two waves originating from the discontinuity in the initial datum suffice - this leads to a single ``turbulent" region with constant density (one notable exception is \cite{MR3831840}). For the particular case of the contact discontinuity with constant pressure $\rho_+=\rho_-$ in Theorem \ref{main_theorem_nonunique} two waves cannot work, even though one might expect parallels to the incompressible vortex sheet case \cite{MR2842999} - this has been noted by several authors. Here we give a short argument to convince the reader why this is the case: 

First of all, note that if $\rho_+=\rho_-$ and in the turbulent zone the pressure is given by a single constant $\rho_1$, then by conservation of mass necessarily $\rho_1=\rho_+=\rho_-$. Therefore the corresponding subsolution $(\rho,v,R)$ has to have constant pressure. However, looking at \eqref{subsolution_system} it is easy to see that then $R\equiv 0$ and the subsolution has to agree to the Riemann solution, ruling out the possibility of convex integration. 

\begin{figure}[tb]
      \includegraphics[width=.99\textwidth]{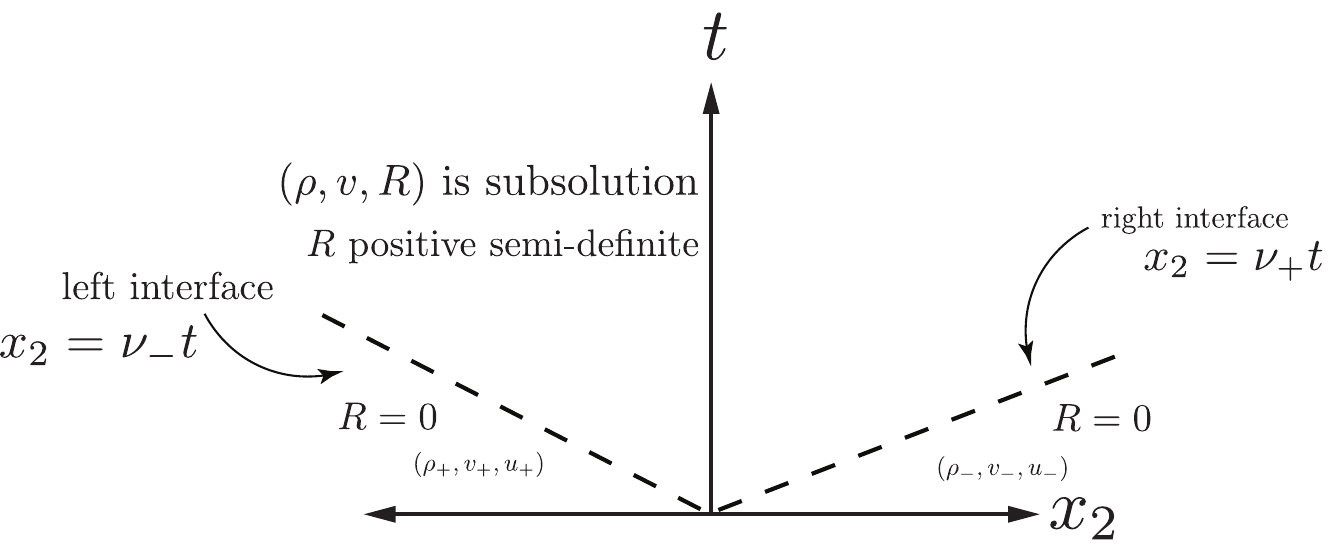}\hspace{.8in}
  \caption{A diagram of a potential 1-D subsolution in space-time (see \eqref{subsolution_system}). In between the left and right interfaces, $(\rho, v,R)$ is a subsolution and $R$ must be positive semi-definite. We also have $\rho=\rho_+=\rho_-=\mathrm{const}$.}\label{fig:subsol}
\end{figure}

More precisely, consider the case $\rho\equiv\mathrm{constant} >0$ and $v$ and $R$ are only functions of $(x_2,t)$, i.e. $v=v(x_2,t)=(v_1(x_2,t),v_2(x_2,t))$ and $R=R(x_2,t)$. From the continuity equation in \eqref{subsolution_system} we deduce $\partial_{x_2} v_2=0$ and thus $v_2$ depends only on $t$, $v_2=v_2(t)$. Then, the second component of the relaxed momentum conservation (second line in \eqref{subsolution_system}) reads
\begin{align}\label{deriv_this}
    \partial_t v_2+\partial_{x_2}(v_2^2+R_{22})=0,
\end{align}
where $R_{22}$ is the component of $R$ in the $(2,2)$ entry. In particular we deduce $\partial_{x_2}^2 R_{22}=0$, which implies that $R_{22}(\cdot,t)$ is a linear function of $x_2$ for each fixed $t$. However, as $x_2\to \pm \infty$ we must have $R\to 0$ (compare with \Cref{fig:subsol}), thus implying that $R_{22}\equiv 0$. But then, being positive semidefinite, we must have $R\equiv 0$.

In fact, in our proof of Theorem \ref{main_theorem_nonunique} we use fan subsolutions consisting of four distinct waves (amounting to 5 different regions where the density is constant) originating from the initial discontinuity. A plot of the waves in our numerical example is depicted in Figure \ref{fig:example_fan}, and the corresponding values of the density are contained in Appendix \ref{sec:numerics}.

\subsubsection{Computational experiments}
Our approach, being computational, is very general and flexible, and our preliminary experimenting shows that we can handle \emph{essentially any class of Riemann initial data.} By using our fan subsolution for the contact discontinuity as a starting point, a numerical search\footnote{Our code, including this numerical experiment, are available on the GitHub at: \url{https://github.com/sammykrupa/NonUniqueness2DIsentropicEuler}} shows numerical solutions to the algebraic system defining a fan subsolution, for Riemann initial data chosen at random: the left-hand state is taken from the box $(\rho_\pm -\frac{1}{2},\rho_\pm+\frac{1}{2})\times(v_{-1}-\frac{1}{2},v_{-1}+\frac{1}{2})\times(v_{\pm2}-\frac{1}{2},v_{\pm2}+\frac{1}{2})$ and the right-hand state is taken from the box $(\rho_\pm -\frac{1}{2},\rho_\pm+\frac{1}{2})\times(v_{+1}-\frac{1}{2},v_{+1}+\frac{1}{2})\times(v_{\pm2}-\frac{1}{2},v_{\pm2}+\frac{1}{2})$. Here, the $(\rho_\pm,v_{\pm1},v_{\pm2})$ denote the Riemann data from \Cref{main_theorem_nonunique}, above. We sampled 100 such Riemann problems at random, and for \emph{all of them} we found a numerical solution of the fan subsolution system\footnote{See the code on GitHub for more details.}. Compare this result with Chen-Chen \cite{MR2303477}, where it is shown that for a certain class of pressure laws, 1-D rarefactions are unique and stable in $L^\infty$, even in the class of truly two-dimensional  solutions to \eqref{system}-\eqref{system_adm} (see also \cite{MR3401974,MR3357629}). We remark that in particular, the work by Chen-Chen \cite{MR2303477} can handle the polytropic pressure law with any $\gamma>1$. As discussed above, for these pressure laws, and for Riemann data whose classical solution contains a shock, non-uniqueness via convex integration is known \cite{MR3744380}.

We thus are led to the conjecture

\begin{conjecture}
Given \emph{any} Riemann initial data \eqref{init_data}, there exists a smooth pressure law $p$, verifying $p'>0$, such that for the system \eqref{system}, \eqref{system_adm} with this pressure law, the Riemann initial data admits infinitely many bounded, admissible solutions.
\end{conjecture}

\vspace{.1in}

\subsubsection{\Cref{main_theorem_unique}: Uniqueness}\label{comments:uniqueness}

We use the weak/strong stability theory of Dafermos and DiPerna (in particular, following \cite[p.~125]{dafermos_big_book}), which is notable for not having any regularity assumptions. The difficulty is that the weak/strong theory cannot allow for jump discontinuities in the solution we want to show uniqueness for. Our proof of \Cref{main_theorem_unique} relies on the fact that under a change of coordinates (Eulerian to Lagrangian), we can ``linearize'' the contact discontinuity of \eqref{system} in some sense and thus a mollification argument will work. Interestingly, this change of coordinates can only be executed in the one-dimensional context.




\subsection{Plan for the paper}
The paper is organized as follows. In \Cref{sec:convex_int}, we give a bare-bones introduction to the convex integration framework we will use.  In \Cref{proof:main_theorem_nonunique}, we give the proof of \Cref{main_theorem_nonunique}. In \Cref{proof:main_theorem_unique}, we prove \Cref{main_theorem_unique}. 

\section{The engine of convex integration}\label{sec:convex_int}
In this section we recall the framework for convex integration for the 2-D isentropic Euler equations as developed by Chiodaroli-De Lellis-Kreml \cite{multi_d_illposed}. We will use $\mathcal{S}_0^{2 \times 2}$ to denote the set of symmetric traceless $2 \times 2$ matrices, and $\mathrm{Id}$ denotes the identity matrix.

\begin{definition}[\protect{Fan Partition \cite[Definition 3.3]{multi_d_illposed}}]\label{def:fan}
A \emph{fan partition} of $\mathbb{R}^2 \times [0, \infty)$ consists of finitely many open sets $P_-, P_1, \ldots, P_N, P_+$ of the following forms:
\begin{align}
P_- &= \{ (x, t) : t > 0 \text{ and } x_2 < \nu_- t \},  \\
P_+ &= \{ (x, t) : t > 0 \text{ and } x_2 > \nu_+ t \},  \\
P_i &= \{ (x, t) : t > 0 \text{ and } \nu_{i-1} t < x_2 < \nu_i t \}, 
\end{align}
where 
\begin{align}\label{fan_ordering}
    \nu_- = \nu_0 < \nu_1 < \cdots < \nu_N = \nu_+
\end{align} 
can be any real numbers.
\end{definition}

\begin{figure}[tb]
      \includegraphics[width=.9\textwidth]{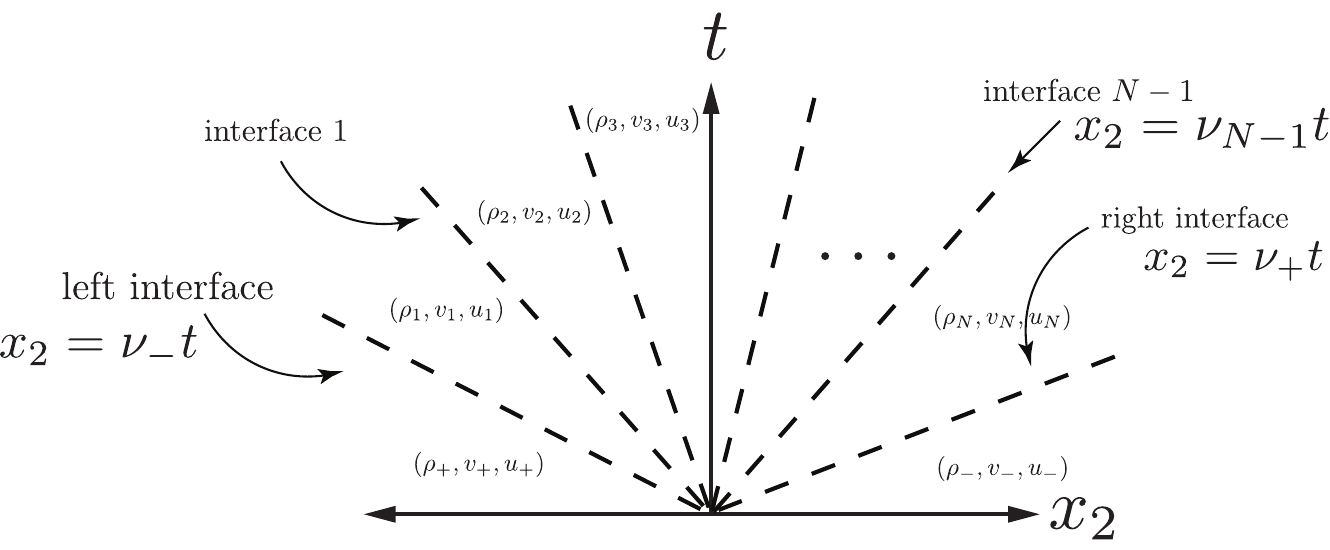}\hspace{.8in}
  \caption{A diagram of a fan partition in space-time.}\label{fig:fan}
\end{figure}

\begin{definition}[\protect{Fan Subsolutions \cite[Definition 3.4]{multi_d_illposed}}]\label{fan_def}
A \emph{fan subsolution} to the system \eqref{system} with initial data \eqref{init_data} is a triple
\[
(\bar{\rho}, \bar{v}, \bar{u}) : \mathbb{R}^2 \times [0, \infty) \to (\mathbb{R}^+, \mathbb{R}^2, \mathcal{S}_0^{2 \times 2})
\]
of piecewise constant functions verifying the following requirements:
\begin{enumerate}
    \item[(i)] There is a fan partition $P_-, P_1, \ldots, P_N, P_+$ of $\mathbb{R}^2 \times [0, \infty)$ such that
    \[
    (\bar{\rho}, \bar{v}, \bar{u}) = \sum_{i=1}^N (\rho_i, v_i, u_i) \mathbbm{1}_{P_i} + (\rho_-, v_-, u_-) \mathbbm{1}_{P_-} + (\rho_+, v_+, u_+) \mathbbm{1}_{P_+}
    \]
    where $\rho_i, v_i, u_i$ are constants with $\rho_i > 0$, $u_{\pm} = v_{\pm} \otimes v_{\pm} - \frac{1}{2} |v_{\pm}|^2 \text{Id}$, and $\rho_{\pm} > 0$; notice that, as $t \to 0^+$, $(\rho(\cdot, t), v(\cdot, t))$ converges to the pair $(\rho^0, v^0)$ of \eqref{init_data}, and this is the reason we say that the triple $(\bar{\rho}, \bar{v}, \bar{u})$ has initial data $(\rho^0, v^0)$.
    \item[(ii)] For every $i \in \{1, 2, \ldots, N\}$ there exists a positive constant $C_i$ such that
    \[
    v_i \otimes v_i - u_i < \frac{C_i}{2} \text{Id}. 
    \]
    \item[(iii)] The triple $(\bar{\rho}, \bar{v}, \bar{u})$ solves the following system in the sense of distributions:
    \begin{align}
    &\partial_t \bar{\rho} + \text{div}_x (\bar{\rho} \bar{v}) = 0, \label{sub_eq1}\\
    &\partial_t (\bar{\rho} \bar{v}) + \text{div}_x (\bar{\rho} \bar{u}) + \nabla_x \left( p(\bar{\rho}) + \frac{1}{2} \left( \sum_i C_i \rho_i \mathbbm{1}_{P_i} + \bar{\rho} |\bar{v}|^2 \mathbbm{1}_{P_- \cup P_+} \right) \right) = 0.\label{sub_eq2}
    \end{align}
\end{enumerate}
\end{definition}

See \Cref{fig:fan}.

\begin{definition}[\protect{Admissible Fan Subsolutions \cite[Definition 3.5]{multi_d_illposed}}]\label{def_fan_adm}

\hfill

A fan subsolution $(\bar{\rho}, \bar{v}, \bar{u})$ is said to be \emph{admissible} if it verifies the following inequality in the sense of distributions
\begin{equation}\label{sub_eq3}
\begin{split}	
\partial_t (\bar{\rho} e(\bar{\rho})) + \text{div}_x &\left[ (\bar{\rho} e(\bar{\rho}) + p(\bar{\rho})) \bar{v} \right] + \partial_t \left( \frac{\bar{\rho} |\bar{v}|^2}{2} \mathbbm{1}_{P_+ \cup P_-} \right) + \text{div}_x \left( \bar{\rho} \frac{|\bar{v}|^2}{2} \bar{v} \mathbbm{1}_{P_+ \cup P_-} \right)\\
&+ \sum_{i=1}^N \left[ \partial_t \left( \rho_i \frac{C_i}{2} \mathbbm{1}_{P_i} \right) + \text{div}_x \left( \rho_i \bar{v} \frac{C_i}{2} \mathbbm{1}_{P_i} \right) \right] \leq 0.
\end{split}
\end{equation}
\end{definition}

\begin{remark}
Fan subsolutions are particular examples of subsolutions as defined in the introduction (cf.~\cite{MR2564474,MR2917063}). Indeed, given a fan subsolution $(\bar{\rho},\bar{v},\bar{u})$ with fan partition $P_-,P_1,\dots P_N,P_+$ and constants $C_i$ as in Definitions \ref{def:fan} and \ref{fan_def}, set
\begin{equation}\label{e:defR}
{R}=\bar{u}-\bar{v}\otimes\bar{v}+\frac{1}{2}\left(\sum_{i=1}^NC_i\mathbbm{1}_{P_i}+|v_-|^2\mathbbm{1}_{P_-}+|v_+|^2\mathbbm{1}_{P_+}\right)\mathrm{Id}.
\end{equation}
Observe that  
\begin{equation*}
\begin{split}
\textrm{tr }{R}&=	\left(\sum_{i=1}^NC_i\mathbbm{1}_{P_i}+|v_-|^2\mathbbm{1}_{P_-}+|v_+|^2\mathbbm{1}_{P_+}\right)-|\bar{v}|^2\\
&=\sum_{i=1}^N(C_i-|v_i|^2)\mathbbm{1}_{P_i},
\end{split}
\end{equation*}
so that $\tfrac{1}{2}(C_i-|v_i|^2)$ can be interpreted as the kinetic energy contained in the high-frequency ``turbulent" part of the flow that is constructed with convex integration. Moreover, conditions (i) and (ii) in Definition \ref{fan_def} imply that $R=0$ in $P_-\cup P_+$ and with $R>0$ (positive definite) in $\bigcup_iP_i$.  
Finally, \eqref{sub_eq1}-\eqref{sub_eq3} follows directly from \eqref{subsolution_system}.
\end{remark}

The main application of admissible fan subsolutions is the following Proposition, which reduces the question of nonunique solutions to \eqref{system}, \eqref{system_adm} to the existence of an admissible fan subsolution.

\begin{proposition}[\protect{Reduction to Admissible Fan Subsolutions \cite[Proposition 3.6]{multi_d_illposed}}]\label{fan_is_enough}

\hfill

Let $p$ be any $C^1$ function and $(\rho_{\pm}, v_{\pm})$ be such that there exists at least one admissible fan subsolution $(\bar{\rho}, \bar{v}, \bar{u})$ of \eqref{system} with initial data \eqref{init_data}. Then there are infinitely many bounded admissible solutions $(\rho, v)$ to \eqref{system}, \eqref{system_adm}, \eqref{init_data} such that $\rho = \bar{\rho}$.
\end{proposition}

In this paper we will consider subsolutions with a fan partition consisting of  5 sets: $P_-$, $P_1,P_2,P_3$ and $P_+$ (see \Cref{fig:example_fan}).

Because the fan subsolution is inherently composed of piecewise constant functions, the PDEs \eqref{sub_eq1} and \eqref{sub_eq2}, along with \Cref{def_fan_adm}, can be equivalently stated as a set of Rankine-Hugoniot conditions at each of the interfaces within the fan partition. More precisely, we have

\begin{proposition}[\protect{\cite[Proposition 5.1]{multi_d_illposed} and \cite[p.~6]{MR3831840}}]\label{RH_prop}
Fix $N\in\mathbb{N}$, $N>0$, and let $P_-, P_1,\ldots,P_N, P_+$ be a fan partition as in \Cref{fan_def}. Consider the constants 
\begin{align*}
v_-, v_+, \rho_-, \rho_+, \rho_1,\ldots,\rho_N,v_1,\ldots,v_N, u_1,\ldots,u_N,C_1,\ldots,C_N
\end{align*}
as in \eqref{3.3}–\eqref{3.14}, where we write  \[
v_i = (\alpha_i, \beta_i),
\]
\[
v_- = (v_{-1}, v_{-2}),
\]
\[
v_+ = (v_{+1}, v_{+2}),
\]
\[
u_i = \begin{pmatrix} \gamma_i & \delta_i \\ \delta_i & -\gamma_i \end{pmatrix},
\]
for $i = 1,\ldots,N$ and some constants $\alpha_i,\beta_i,\gamma_i,\delta_i,v_{-1}, v_{-2},v_{+1}, v_{+2}$.

Then, we have an admissible fan subsolution (as in \Cref{fan_def} and \Cref{def_fan_adm}) if and only if the following equalities and inequalities hold:


\begin{itemize}

 \item Rankine-Hugoniot conditions on the left interface:
    \begin{align}
    \nu_- (\rho_- - \rho_1) &= \rho_- v_{-2} - \rho_1 \beta_1, \label{3.1} \\
    \nu_- (\rho_- v_{-1} - \rho_1 \alpha_1) &= \rho_- v_{-1} v_{-2} - \rho_1 \delta_1, \label{3.2} \\
    \nu_- (\rho_- v_{-2} - \rho_1 \beta_1) &= \rho_- v_{-2}^2 + \rho_1 \gamma_1 + p(\rho_-) - p(\rho_1) - \rho_1 \frac{C_1}{2}; \label{3.3}
    \end{align}
   
    \item Rankine-Hugoniot conditions on interface $i$ ($i \in \{1, \ldots, N - 1\}$):
    \begin{align}
    \nu_i (\rho_i - \rho_{i+1}) &= \rho_i \beta_i - \rho_{i+1} \beta_{i+1}, \label{3.4} \\
    \nu_i (\rho_i \alpha_i - \rho_{i+1} \alpha_{i+1}) &= \rho_i \delta_i - \rho_{i+1} \delta_{i+1}, \label{3.5} \\
    \nu_i (\rho_i \beta_i - \rho_{i+1} \beta_{i+1}) &= -\rho_i \gamma_i + \rho_{i+1} \gamma_{i+1} + p(\rho_i) - p(\rho_{i+1}) + \rho_i \frac{C_i}{2} - \rho_{i+1} \frac{C_{i+1}}{2}; \label{3.6}
    \end{align}
    \item Rankine-Hugoniot conditions on the right interface:
    \begin{align}
    \nu_+ (\rho_N - \rho_+) &= \rho_N \beta_N - \rho_+ v_{+2}, \label{3.7} \\
    \nu_+ (\rho_N \alpha_N - \rho_+ v_{+1}) &= \rho_N \delta_N - \rho_+ v_{+1} v_{+2}, \label{3.8} \\
    \nu_+ (\rho_N \beta_N - \rho_+ v_{+2}) &= -\rho_N \gamma_N - \rho_+ v_{+2}^2 + p(\rho_N) - p(\rho_{+}) + \rho_N \frac{C_N}{2}; \label{3.9}
    \end{align}
    
    \item Subsolution conditions ($i \in \{1, \ldots, N\}$):
    \begin{align}
    &\alpha_i^2 + \beta_i^2 < C_i \label{3.10} \\
    &\left( \frac{C_i}{2} - \alpha_i^2 + \gamma_i \right) \left( \frac{C_i}{2} - \beta_i^2 - \gamma_i \right) - (\delta_i - \alpha_i \beta_i)^2 > 0 \label{3.11}
    \end{align}
    \item Admissibility condition on the left interface:
    \begin{multline}\label{3.12}
    \nu_- (\rho_- \epsilon(\rho_-) - \rho_1 \epsilon(\rho_1)) + \nu_- \left( \rho_- \frac{|v_-|^2}{2} - \rho_1 \frac{C_1}{2} \right)  
    \\
   \leq [ (\rho_- \epsilon(\rho_-) + p(\rho_-)) v_{-2} - (\rho_1 \epsilon(\rho_1) + p(\rho_1)) \beta_1 ] 
    \\
    \quad + \left( \rho_- v_{-2} \frac{|v_-|^2}{2} - \rho_1 \beta_1 \frac{C_1}{2} \right); 
    \end{multline}
    \item Admissibility condition on interface $i$ ($i \in \{1, \ldots, N - 1\}$):
    
    \begin{multline}\label{3.13}
    \nu_i (\rho_i \epsilon(\rho_i) - \rho_{i+1} \epsilon(\rho_{i+1})) + \nu_i \left( \rho_i \frac{C_i}{2} - \rho_{i+1} \frac{C_{i+1}}{2} \right) 
    \\
   \leq  [ (\rho_i \epsilon(\rho_i) + p(\rho_i)) \beta_i - (\rho_{i+1} \epsilon(\rho_{i+1}) + p(\rho_{i+1})) \beta_{i+1} ] 
    \\
    \quad + \left( \rho_i \beta_i \frac{C_i}{2} - \rho_{i+1} \beta_{i+1} \frac{C_{i+1}}{2} \right); 
    \end{multline}

    \item Admissibility condition on the right interface:
    \begin{multline}\label{3.14}
    \nu_+ (\rho_N \epsilon(\rho_N) - \rho_+ \epsilon(\rho_+)) + \nu_+ \left( \rho_N \frac{C_N}{2} - \rho_+ \frac{|v_+|^2}{2} \right) 
    \\
    \leq [ (\rho_N \epsilon(\rho_N) + p(\rho_N)) \beta_N - (\rho_+ \epsilon(\rho_+) + p(\rho_+)) v_{+2} ] 
    \\
    \quad + \left( \rho_N \beta_N \frac{C_N}{2} - \rho_+ v_{+2} \frac{|v_+|^2}{2} \right). 
    \end{multline}
\end{itemize}
\end{proposition}

\section{Proof of \Cref{main_theorem_nonunique}}\label{proof:main_theorem_nonunique}

We first state our key ``MATLAB'' Lemma, which states the result of our symbolic MATLAB code. It is not immediately clear why we need some of the algebraic relations which we state. They will become clear once we begin the proof of \Cref{main_theorem_nonunique}.

\begin{lemma}[The MATLAB Lemma]\label{MATLAB_prop}
    With $N=3$, we can find constants 
    \begin{multline}
    \hat v_{-1}, \hat v_{-2}, \hat v_{+1}, \hat v_{+2}, \hat \rho_-, \hat \rho_+, \hat \rho_1,\ldots,\hat \rho_N,\hat \nu_-,\hat \nu_+,\hat \nu_1,\ldots,\hat \nu_{N-1}, \hat \alpha_1,\ldots,\hat \alpha_N,\hat \beta_1,\ldots,\hat \beta_N,
    \\
\hat \gamma_1,\ldots,\hat \gamma_N,\hat \delta_1,\ldots,\hat \delta_N,\hat C_1,\ldots,\hat C_N
    \end{multline}
which, playing the role of 
\begin{multline}
     v_{-1},  v_{-2},  v_{+1},  v_{+2},  \rho_-,  \rho_+,  \rho_1,\ldots, \rho_N, \nu_-, \nu_+, \nu_1,\ldots, \nu_{N-1},  \alpha_1,\ldots, \alpha_N, \beta_1,\ldots, \beta_N,
    \\
 \gamma_1,\ldots, \gamma_N, \delta_1,\ldots, \delta_N, C_1,\ldots, C_N,
    \end{multline}
    respectively, will satisfy \eqref{3.1}, \eqref{3.2}, \eqref{3.3}, as well as \eqref{3.10}-\eqref{3.14}, \eqref{fan_ordering}, and \eqref{3.4}-\eqref{3.6} for $i=1$.

Furthermore, \eqref{3.4}-\eqref{3.6} for $i=2$ and \eqref{3.7}-\eqref{3.9} are satisfied approximately: for each equation, the magnitude of the difference between the right-hand side and the left-hand side is always less than $10^{-11}$.

In equations \eqref{3.12}-\eqref{3.14}, the terms $\epsilon(\rho_-),\epsilon(\rho_+),\epsilon(\rho_1),\ldots,\epsilon(\rho_N)$ are substituted, respectively, by constants $\hat \epsilon_-,\hat\epsilon_+,\hat\epsilon_1,\ldots,\hat\epsilon_N$.

Similarly, in equations \eqref{3.1}-\eqref{3.14}, the terms $p(\rho_-),p(\rho_+),p(\rho_1),\ldots,p(\rho_N)$ are substituted by the quantities $(\hat \rho_-)^2 \hat \epsilon_-',(\hat \rho_+)^2 \hat \epsilon_+',(\hat \rho_1)^2 \hat \epsilon_1',\ldots,(\hat \rho_N)^2 \hat \epsilon_N'$, respectively, for constants $\hat \epsilon_-', \hat \epsilon_+',\hat \epsilon_1',\hat \epsilon_N'$ (cf. \eqref{pressure_relation}).

For example, \eqref{3.14} becomes

\begin{multline*}
        \hat\nu_+ (\hat\rho_N \hat\epsilon_N - \hat\rho_+ \hat\epsilon_+) + \hat\nu_+ \left( \hat\rho_N \frac{\hat C_N}{2} - \hat\rho_+ \frac{|(\hat v_{+1},\hat v_{+2})|^2}{2} \right) \leq 
        \\
        [ (\hat\rho_N \hat\epsilon_N + (\hat \rho_N)^2\hat\epsilon_N') \hat\beta_N - (\hat \rho_+ \hat \epsilon_+ + (\hat \rho_+)^2\hat\epsilon_+')) \hat v_{+2} ]  \\
    \quad + \left( \hat\rho_N \hat\beta_N \frac{\hat C_N}{2} - \hat \rho_+ \hat v_{+2} \frac{|(\hat v_{+1},\hat v_{+2})|^2}{2} \right).
\end{multline*}

We also have the following equalities
\begin{equation*}
    \hat\rho_+=\hat\rho_-,\quad \hat\epsilon_+=\hat\epsilon_-,\quad \hat\epsilon_+'=\hat\epsilon_-',\quad \hat v_{+2}=\hat v_{-2},
\end{equation*}
while
\begin{equation*}
    \hat v_{+1}\neq\hat v_{-1}.
\end{equation*}

Moreover, we have the additional set of inequalities

\begin{align}\label{convex_ineq}
\hat\epsilon_j-\hat\epsilon_i-\hat\epsilon_i'(\hat \rho_j-\hat\rho_i)>0,
\end{align}
for all $i,j\in\{0,1,\ldots,N\}$ ($i\neq j$), where $\hat\rho_0\coloneqq \hat\rho_-$, $\hat\epsilon_0\coloneqq \hat\epsilon_-$, and $\hat\epsilon_0'\coloneqq \hat\epsilon_-'$. 

Furthermore, 
    \begin{equation} 
    \begin{aligned}\label{positive_stuff}
        \hat\rho_-,\hat\rho_+,  \hat\rho_1,\ldots, \hat\rho_N>0,\\
        \hat \epsilon_-', \hat \epsilon_+',\hat \epsilon_1',\hat \epsilon_N' >0.
    \end{aligned}
    \end{equation}
    
    Every inequality is verified strictly, with the difference between the left-hand side and the right-hand side of the inequality being at least $\frac{1}{3}$. This includes \eqref{convex_ineq}, \eqref{positive_stuff}, \eqref{fan_ordering}, and \eqref{3.10}-\eqref{3.14}.

    Lastly, we have the following bounds:
    \begin{equation}
    \begin{aligned}\label{numerical_bounds}
        \max\{\hat v_{-1}, \hat v_{-2}, \hat v_{+1}, \hat v_{+2}\}&\leq 59\\
        \max\{\abs{\hat \alpha_1},\ldots,\abs{\hat \alpha_N}\}&\leq 59 \\
        \max\{\abs{\hat \beta_1},\ldots,\abs{\hat \beta_N}\}&\leq 19 \\
       \max\{\abs{\hat \gamma_1},\ldots,\abs{\hat \gamma_N}\}&\leq 1531 \\
    \max\{\abs{\hat \delta_1},\ldots,\abs{\hat \delta_N}\}&\leq 1046 \\
     \max\{\abs{\hat\nu_-},\abs{\hat \nu_+},\abs{\hat \nu_1},\ldots,\abs{\hat \nu_{N-1}}\}&\leq 34 \\
     \max\{\abs{\hat C_1},\ldots,\abs{\hat C_N}\}&\leq 3714 \\
      \max\{\abs{\hat\rho_-},\abs{\hat \rho_+},\abs{\hat \rho_1},\ldots,\abs{\hat \rho_{N}}\}&\leq 13 \\
      \max\{\abs{\hat\epsilon_-},\abs{\hat \epsilon_+},\abs{\hat \epsilon_1},\ldots,\abs{\hat \epsilon_{N}}\}&\leq 2308 \\
    \max\{\abs{\hat\epsilon_-'},\abs{\hat \epsilon_+'},\abs{\hat \epsilon_1'},\ldots,\abs{\hat \epsilon_{N}'}\}&\leq 7\\
    \sigma_{\text{min}}(D\Gamma(\hat\alpha_N,\hat\beta_N,\hat\delta_N,\hat\rho_N,\hat\nu_+,\hat\nu_{N-1})) &\geq 2.
    \end{aligned}
     \end{equation}
     where $\sigma_{\text{min}}(D\Gamma)$ denotes the smallest singular value of the matrix $D\Gamma$. The matrix $D\Gamma$ (see \eqref{DGamma}, below) is the Jacobian of the function $\Gamma$, defined below (see \eqref{Gamma}).

\end{lemma}

The proof of \Cref{MATLAB_prop} is in \Cref{proof:matlab}, below. We give explicit constants which verify the conclusions of \Cref{MATLAB_prop} in the Appendix (see \Cref{sec:numerics}). 

\vspace{.05in}

We can now start the proof of \Cref{main_theorem_nonunique}.

\uline{Step 1} \Cref{MATLAB_prop} gives approximate solutions to \eqref{3.4}-\eqref{3.6} for $i=2$ and \eqref{3.7}-\eqref{3.9}: for each equation, the magnitude of the difference between the right-hand side and the left-hand side is always less than $10^{-11}$. Let us first perturb this approximate solution into an exact solution, using the Inverse Function Theorem. In particular, we will use a quantitative version of the Inverse Function Theorem (see \Cref{sec:IFT}).

Consider the map $\Gamma\colon\mathbb{R}^6\to\mathbb{R}^6$, given by
\begin{multline}\label{Gamma}
\Gamma\colon (\alpha,\beta,\delta,\rho,\nu,\tilde{\nu})\mapsto
    \\
    \begin{bmatrix}
       \hat\rho_{N-1}\hat\beta_{N-1} -\rho\beta-\tilde{\nu}(\hat\rho_{N-1}-\rho)\\[0.3em]
        \hat\rho_{N-1}\hat\delta_{N-1}-\rho\delta-\tilde{\nu}(\hat\rho_{N-1}\hat\alpha_{N-1}-\rho\alpha)\\[0.3em]
       -\hat\rho_{N-1}\hat\gamma_{N-1}+\rho\hat\gamma_{N}+(\hat\rho_{N-1})^2\hat\epsilon_{N-1}'-\rho^2 \hat\epsilon_{N}'+\hat\rho_{N-1}\frac{\hat C_{N-1}}{2}-\rho\frac{\hat C_{N}}{2}-\tilde{\nu}(\hat\rho_{N-1}\hat\beta_{N-1}-\rho\beta)\\[0.3em]
       \rho\beta-\hat\rho_{+}\hat v_{+2}-\nu(\rho-\hat\rho_{+})\\[0.3em]
       \rho\delta-\hat\rho_{+}\hat v_{+1}\hat v_{+2}-\nu(\rho\alpha -\hat\rho_{+}\hat v_{+1})\\[0.3em]
       -\rho\hat\gamma_{N}-\hat\rho_{+}(\hat v_{+2})^2+\rho^2\hat\epsilon_{N}'-(\hat\rho_{+})^2\hat\epsilon_{+}'+\rho\frac{\hat C_N}{2}-\nu (\rho\beta-\hat\rho_+\hat v_{+2})
     \end{bmatrix}.
\end{multline}

The first three rows of $\Gamma$ correspond with \eqref{3.4}-\eqref{3.6} (with $i=2$) and the last three rows of $\Gamma$ correspond with \eqref{3.7}-\eqref{3.9}.

We also calculate the Jacobian of $\Gamma$:

\begin{multline}\label{DGamma}
    D\Gamma(\alpha,\beta,\delta,\rho,\nu,\tilde{\nu}) = 
    \\
     \begin{bmatrix}
       0 & -\rho & 0   &  -\beta+\tilde{\nu} &  0  &   -(\hat\rho_{N-1}-\rho)       \\[0.3em]
       \rho\tilde{\nu} & 0           & -\rho  & -\delta +\alpha\tilde{\nu} & 0  & -(\hat\rho_{N-1}\hat\alpha_{N-1}-\rho\alpha) \\[0.3em]
       0           & \tilde{\nu}\rho & 0    & \hat\gamma_{N}-2\rho\hat\epsilon_N'-\frac{\hat C_N}{2}+\beta\tilde{\nu} & 0 & -(\hat\rho_{N-1}\hat\beta_{N-1}-\rho\beta) \\[0.3em]
       0 & \rho & 0  & \beta-\nu & -(\rho-\hat\rho_+) & 0         \\[0.3em]
       -\nu\rho & 0           & \rho & \delta-\nu\alpha & -(\rho\alpha-\hat\rho_+\hat v_{+1}) & 0 \\[0.3em]
       0           & -\nu\rho & 0 & -\hat\gamma_N+2\rho\hat\epsilon_N'+\frac{\hat C_N}{2}-\nu\beta & -(\rho\beta-\hat\rho_+\hat v_{+2}) & 0
       \end{bmatrix}.
\end{multline}

Using \Cref{IFT_prop} with the bounds \eqref{numerical_bounds} we can perturb the values of 
\begin{equation*}
	\hat\alpha_N,\hat\beta_N,\hat\delta_N,\hat\rho_N,\hat\nu_{+},\hat\nu_{N-1},
\end{equation*}
 keeping all other constants fixed, such that \eqref{3.4}-\eqref{3.6} (for $i=2$) and \eqref{3.7}-\eqref{3.9} are verified exactly. 


Keeping the other constants 
    \begin{multline}
    \hat v_{-1}, \hat v_{-2}, \hat v_{+1}, \hat v_{+2}, \hat \rho_-, \hat \rho_+, \hat \rho_1,\ldots,\hat \rho_{N-1},\hat \nu_-,\hat \nu_1,\ldots,\hat \nu_{N-2}, \hat \alpha_1,\ldots,\hat \alpha_{N-1},\hat \beta_1,\ldots,\hat \beta_{N-1},
    \\
\hat \gamma_1,\ldots,\hat \gamma_N,\hat \delta_1,\ldots,\hat \delta_{N-1},\hat C_1,\ldots,\hat C_N
    \end{multline}
fixed, this gives us exact values of
\begin{multline}
     v_{-1},  v_{-2},  v_{+1},  v_{+2},  \rho_-,  \rho_+,  \rho_1,\ldots, \rho_N, \nu_-, \nu_+, \nu_1,\ldots, \nu_{N-1},  \alpha_1,\ldots, \alpha_N, \beta_1,\ldots, \beta_N,
    \\
 \gamma_1,\ldots, \gamma_N, \delta_1,\ldots, \delta_N, C_1,\ldots, C_N,
    \end{multline}
such that we verify \eqref{3.1}-\eqref{3.6} for $i=1,2$, and \eqref{3.7}-\eqref{3.9}. By making the $A$ (in the context of \Cref{IFT_prop}) large enough, and also from \eqref{numerical_bounds}, we verify \eqref{convex_ineq}, \eqref{positive_stuff}, \eqref{fan_ordering}, and \eqref{3.10}-\eqref{3.14}. Here, we are assuming the existence of functions $\epsilon$ and $p$ such that $\epsilon(\rho_i)=\hat\epsilon_i$, $\epsilon'(\rho_i)=\hat\epsilon_i'$  for all $i$ and $\epsilon(\rho_\pm)=\hat\epsilon_\pm$, $\epsilon'(\rho_\pm)=\hat\epsilon_\pm'$ and $p(\rho)=\rho^2\epsilon'(\rho)$ for all $\rho\in(0,\infty)$. We will construct these functions next.

\uline{Step 2}

We will take advantage of the following standard fact about convex functions (cf. \eqref{convex_ineq}, above):

\begin{lemma}[\protect{Algebraic inequalities which yield a strictly convex function}]\label{convex_extension_prop}

\hfill

Fix $n, N \in \mathbb{N}$. Assume there exists $(x_i, h_i, D_i) \in \mathbb{R}^n \times \mathbb{R} \times \mathbb{R}^n$, for $i = 1, \ldots, N$, such that the strict inequalities
\begin{equation}
h_j > h_i + D_i \cdot (x_j - x_i) \quad \text{for all } i \neq j
\end{equation}
are verified. Then, there exists a smooth and strictly convex function $\xi\colon  \mathbb{R}^n \to \mathbb{R}$ such that $\xi(x_i) = h_i$ and $D\xi(x_i) = D_i$, for all $i$.
\end{lemma}
\begin{remark}
It is well known that such an $\xi$ exists if we only ask for it to be convex. See e.g. \cite[p.~143]{MR2048569}. However, in this \Cref{convex_extension_prop} $\xi$ is constructed in such a way to be \emph{strictly} convex. 
\end{remark}

The proof of \Cref{convex_extension_prop} is standard and can be found in e.g. \cite[p.~16]{2024arXiv240307784K}.

\begin{figure}[tb]
\centering
      \includegraphics[width=\textwidth]{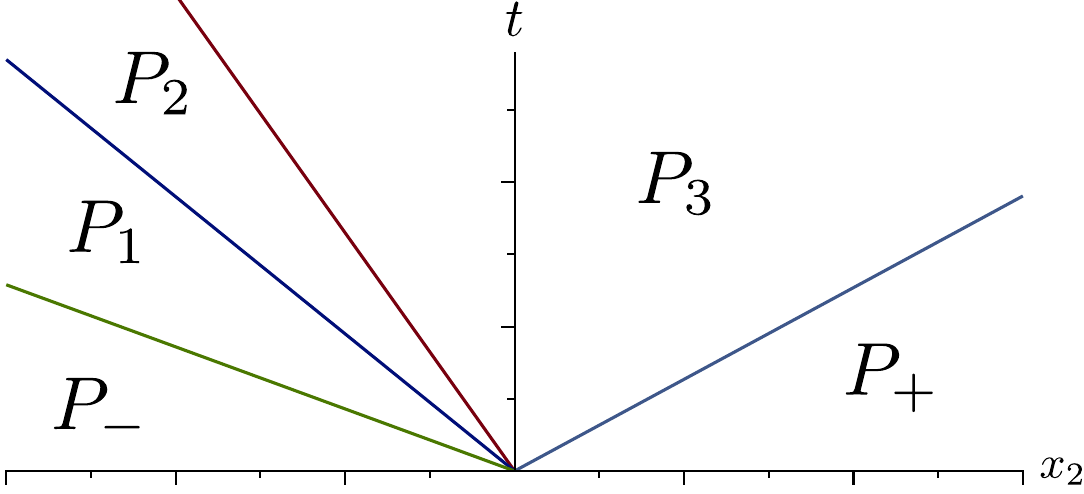}
  \caption{The fan partition used for our fan subsolution with contact discontinuity initial data, with the slopes of the lines determined by \Cref{MATLAB_prop}.}\label{fig:example_fan}
\end{figure}

From \eqref{convex_ineq} and \Cref{convex_extension_prop}, we can find a smooth, strictly convex function $\epsilon:\mathbb{R}\to\mathbb{R}$ such that  $\epsilon(\rho_i)=\hat\epsilon_i$, $\epsilon'(\rho_i)=\hat\epsilon_i'$  for all $i$ and $\epsilon(\rho_\pm)=\hat\epsilon_\pm$, $\epsilon'(\rho_\pm)=\hat\epsilon_\pm'$. We can then define
\begin{align}
    p(\rho)\coloneqq \rho^2 \epsilon'(\rho).
\end{align}

Remark that in \Cref{MATLAB_prop}, in equations \eqref{3.1}-\eqref{3.14}, the terms 
\begin{align}
p(\rho_-),p(\rho_+),p(\rho_1),\ldots,p(\rho_N)
\end{align} 
are substituted by the quantities 
\begin{align}
(\hat \rho_-)^2 \hat \epsilon_-',(\hat \rho_+)^2 \hat \epsilon_+',(\hat \rho_1)^2 \hat \epsilon_1',\ldots,(\hat \rho_N)^2 \hat \epsilon_N',
\end{align}
cf. \eqref{pressure_relation}.
Thus, we have found functions $p,\epsilon$ and constants such that by \Cref{RH_prop}, we have found an admissible fan subsolution (see \Cref{fan_def}, \Cref{def_fan_adm}). Remark also that $p'(\rho)=2\rho\epsilon'(\rho)+\rho^2\epsilon''(\rho)>0$ by \eqref{positive_stuff}, \eqref{pressure_relation} and the convexity of $\epsilon$.
By \Cref{fan_is_enough}, this proves \Cref{main_theorem_nonunique}. See \Cref{fig:example_fan} for a plot showing the fan partition in the $x_2$-$t$ plane.

\subsection{Proof of \Cref{MATLAB_prop}}\label{proof:matlab}
    The proof is computational, using MATLAB \footnote{Our code is available on the GitHub: \url{https://github.com/sammykrupa/NonUniqueness2DIsentropicEuler}}.

    \uline{Step 1} First, we enter into MATLAB all of the equalities and inequalities we want to satisfy: \eqref{3.1}-\eqref{3.14}, \eqref{convex_ineq}, \eqref{positive_stuff}, and \eqref{fan_ordering}.

    \uline{Step 2} Then, we run the MATLAB R2023b solver \emph{fmincon} and the \emph{interior-point} algorithm (see \cite{MATLAB_solve,fmincon}).  This returns approximate \emph{numeric} values in double precision (for details, see \cite{MATLAB_double}), which we then convert to \emph{symbolic values} within MATLAB (see \cite{MATLAB_symbolic}). Let us call these symbolic values the 
    \begin{multline}\label{tilde_values}
         \tilde v_{-1}, \tilde v_{-2}, \tilde v_{+1}, \tilde v_{+2}, \tilde \rho_-, \tilde \rho_+, \tilde \rho_1,\ldots,\tilde \rho_N,\tilde \nu_-,\tilde \nu_+,\tilde \nu_1,\ldots,\tilde \nu_{N-1}, \tilde \alpha_1,\ldots,\tilde \alpha_N,\tilde \beta_1,\ldots,\tilde \beta_N,
    \\
\tilde \gamma_1,\ldots,\tilde \gamma_N,\tilde \delta_1,\ldots,\tilde \delta_N,\tilde C_1,\ldots,\tilde C_N, \tilde \epsilon_-,\tilde\epsilon_+,\tilde\epsilon_1,\ldots,\tilde\epsilon_N, \tilde \epsilon_-', \tilde \epsilon_+',\tilde \epsilon_1',\tilde \epsilon_N'.
    \end{multline}

    Notice that solving a nonlinear system of inequalities with the interior-point algorithm gives different approximate solutions depending on the initial point at which the algorithm starts. Our code chooses a random initial point for the interior-point algorithm each time the code is run. Remark that not all initial points will converge to a feasible solution. We also include in the code the starting point which yields the symbolic values from \Cref{MATLAB_prop} (see \Cref{sec:numerics}). Moreover, our code also includes the symbolic values themselves from \Cref{MATLAB_prop} if the reader wants to examine our particular solution. 

From here on, we perform all computations symbolically within MATLAB (see \cite{MATLAB_symbolic}). This yields rigorous mathematical statements.

\uline{Step 2}

    Because the MATLAB solver only returns approximate solutions to our algebraic system of equalities and inequalties, the equalities \eqref{3.1}-\eqref{3.9} are not satisfied exactly. Roughly speaking, we can make some of these equalities exact by using the equalities themselves and solving the equalities for a particular constant, and then using the values of constants from \eqref{tilde_values} to determine this particular constant. More precisely, this ``correction algorithm'' works as follows:

For the left interface:

From \eqref{3.1}, we define
\begin{align}
    \hat \nu_-\coloneqq\frac{\tilde\rho_- \tilde v_{-2}-\tilde\rho_1\tilde\beta_1}{\tilde\rho_--\tilde\rho_1}.
\end{align}

From \eqref{3.2}, we define
\begin{align}
    \hat\delta_1\coloneqq\frac{\tilde\nu_-(\tilde\rho_-\tilde v_{-1}-\tilde\rho_1\tilde\alpha_1)-\tilde\rho_1\tilde v_{-1}\tilde v_{-2}}{-\tilde\rho_1}.
\end{align}

From \eqref{3.3}, we define
\begin{align}
    \hat\epsilon_1'\coloneqq \frac{-(\tilde\nu_-(\tilde\rho_-\tilde v_{-2}-\tilde\rho_1\tilde\beta_1)+(\tilde\rho_-(\tilde v_{-2})^2+\tilde\rho_1\tilde\gamma_1+(\tilde\rho_-)^2\tilde\epsilon_-'-\tilde\rho_1\frac{\tilde C_1}{2})}{(\tilde\rho_1)^2}.
\end{align}

Then, for each of the interfaces $i=1,\ldots,N-2$, we do the following (which in the case $N=3$, simply means $i=1$):

From \eqref{3.4}, we define
\begin{align}
    \hat\nu_i\coloneqq\frac{\tilde\rho_i\tilde\beta_i-\tilde\rho_{i+1}\tilde\beta_{i+1}}{\tilde\rho_i-\tilde\rho_{i+1}}.
\end{align}

From \eqref{3.5}, we write
\begin{align}
    \hat\delta_{i+1}\coloneqq \frac{\tilde\nu_i(\tilde\rho_i\tilde\alpha_i-\tilde\rho_{i+1}\tilde\alpha_{i+1})-\tilde\rho_i\tilde\delta_i}{-\tilde\rho_{i+1}}.
\end{align}

From \eqref{3.6}, we write
\begin{align}
     \hat\epsilon_{i+1}'\coloneqq \frac{-\tilde\nu_i(\tilde\rho_i\tilde\beta_i-\tilde\rho_{i+1}\tilde\beta_{i+1})+(-\tilde\rho_i\tilde\gamma_i+\tilde\rho_{i+1}\tilde\gamma_{i+1}+(\tilde\rho_i)^2\tilde\epsilon_i'+\tilde\rho_i\frac{\tilde C_i}{2}-\tilde\rho_{i+1}\frac{\tilde C_{i+1}}{2})}{(\tilde\rho_{i+1})^2}.
\end{align}

 For all other constants we use the symbolic (exact) values given in \eqref{tilde_values}. That is, we set
 \begin{align*}
     \hat v_{-1} &\coloneqq \tilde v_{-1},\quad \hat v_{-2} \coloneqq \tilde v_{-2},\quad \hat v_{+1} \coloneqq \tilde v_{+1},\quad \hat v_{+2} \coloneqq \tilde v_{+2}\\
     \hat \rho_- &\coloneqq \tilde \rho_-,\quad \hat \rho_+ \coloneqq \tilde \rho_+,\quad \hat \rho_1 \coloneqq \tilde \rho_1,\quad \hat \rho_2 \coloneqq \tilde \rho_2,\quad \hat \rho_3 \coloneqq \tilde \rho_3\\
     \hat \nu_+ &\coloneqq \tilde \nu_+,\quad \hat \nu_{2} \coloneqq \tilde \nu_{2}\\
     \hat \alpha_1 &\coloneqq \tilde \alpha_1,\quad \hat \alpha_2 \coloneqq \tilde \alpha_2,\quad \hat \alpha_3 \coloneqq \tilde \alpha_3\\
     \hat \beta_1 &\coloneqq \tilde \beta_1,\quad \hat \beta_2 \coloneqq \tilde \beta_2,\quad \hat \beta_3 \coloneqq \tilde \beta_3\\
\hat \gamma_1 &\coloneqq \tilde \gamma_1,\quad \hat \gamma_2 \coloneqq \tilde \gamma_2,\quad \hat \gamma_3 \coloneqq \tilde \gamma_3\\
\hat \delta_3 &\coloneqq \tilde \delta_3 \\
\hat C_1 &\coloneqq \tilde C_1,\quad \hat C_2 \coloneqq \tilde C_2,\quad \hat C_3 \coloneqq \tilde C_3\\
\hat\epsilon_- &\coloneqq \tilde\epsilon_-,\quad \hat\epsilon_+ \coloneqq \tilde\epsilon_+,\quad \hat\epsilon_1 \coloneqq \tilde\epsilon_1,\quad \hat\epsilon_2 \coloneqq \tilde\epsilon_2,\quad \hat\epsilon_3 \coloneqq \tilde\epsilon_3\\
\hat \epsilon_-' &\coloneqq \tilde \epsilon_-',\quad \hat \epsilon_+' \coloneqq \tilde \epsilon_+',\quad \hat \epsilon_3' \coloneqq \tilde \epsilon_3'.
 \end{align*}

 This gives us constants which symbolic verifications in MATLAB show satisfy the conclusions of the Proposition.

\section{Proof of \Cref{main_theorem_unique}}\label{proof:main_theorem_unique}
The solution $(\rho,v_1,v_2)$ will verify the system of conservation laws in one spatial dimension,
\begin{align}
       \begin{cases}\label{1-D_version}
    \partial_t \rho+\partial_{x_2} (v_2\rho)=0,\\
    \partial_t (v_1\rho) +\partial_{x_2} (v_1 v_2\rho)=0,\\
    \partial_t (v_2\rho) +\partial_{x_2} (\rho v_2^2+p(\rho))=0,
    \end{cases}
    \end{align}
    and it will verify the entropy inequality $\partial_t \eta(\rho,v_1,v_2)+\partial_{x_2} q(\rho,v_1,v_2)\leq 0$ (in a weak sense), for the entropy function
    \begin{align}
        \eta(\rho,v_1,v_2)=\rho\epsilon(\rho)+\frac{v_1^2+v_2^2}{2}\rho,
    \end{align}
    and entropy-flux
    \begin{align}
        q(\rho,v_1,v_2)=\big(\rho\epsilon(\rho)+\frac{v_1^2+v_2^2}{2}\rho +p(\rho)\big)v_2.
    \end{align}

\uline{Step 1}

In the theory of conservation laws in one spatial dimension, there is a one-to-one correspondence between the weak, bounded solutions to an equation in Eulerian coordinates, uniformly away from vacuum, and the weak, bounded solutions to the corresponding equation in Lagrangian coordinates (see \cite[Theorem 2]{wagnergasdynamics}). The correspondence is induced by the Lipschitz-continuous map $T\colon (x_2,t)\mapsto (y(x_2,t),t)$ (for $(x_2,t)\in\mathbb{R}\times[0,\infty)$). Applying this result, we switch \eqref{1-D_version} from Eulerian to Lagrangian coordinates, which gives us
\begin{align}
    \begin{cases}\label{1-D_Lagrangian}
    \partial_t (\frac{1}{\rho}) -\partial_{y} v_2=0,\\
    \partial_t v_1 =0,\\
    \partial_t v_2 +\partial_{y} p(\rho)=0,
    \end{cases}
\end{align}
with a corresponding entropy function given by
\begin{align}
    \bar{\eta}=\epsilon(\rho)+\frac{v_1^2+v_2^2}{2},
\end{align}
with associated entropy-flux
\begin{align}
    \bar{q}=p(\rho)v_2.
\end{align}

Moreover, \cite[Theorem 2]{wagnergasdynamics} tells us that the entropy inequality also holds in Lagrangian coordinates, 
\begin{align}\label{entropy_Lagrangian}
    \partial_t \bar{\eta}+\partial_{y} \bar{q}\leq 0,
\end{align}
and $\bar{\eta}$ is convex in the variables $(w,v_1,v_2)$, where $w\coloneqq \frac{1}{\rho}$ is the \emph{specific volume}. Notice that the (strict) convexity of $\epsilon$, as a function of $w$, also follows directly from $p,p'>0$ and \eqref{pressure_relation}.

Remark that by the one-to-one correspondence, uniqueness for \eqref{1-D_Lagrangian} will imply uniqueness for \eqref{1-D_version}.

\uline{Step 2}

We now apply the theory of weak/strong stability to \eqref{1-D_Lagrangian}, following Dafermos \cite[p.~122]{dafermos_big_book}. 

In particular, let us recall the key estimate in the weak/strong theory,

\begin{multline}\label{key_dafermos}
\int_{|y| < r + s(t - \sigma)} \bar{\eta}(U(y, \sigma)| \bar{U}(y, \sigma)) \, dy \leq \int_{|y| < r + st} \bar{\eta}(U(y,0)| \bar{U}(y,0)) \, dy
\\
- \int_0^\sigma \int_{|y| < r + s(t - \tau)}  \partial_{y} \bar{U} \nabla^2 \bar{\eta}(\bar{U})G(U| \bar{U}) \, dy d\tau,
\end{multline}

where $G$ is the flux for the system \eqref{1-D_Lagrangian}, written as $G(U)=(-v_2,0,p(\rho))$, $\eta(\cdot|\cdot)$ is the relative entropy, defined as
 \begin{align}
     \bar{\eta}(a|b)\coloneqq \bar{\eta}(a)-\bar{\eta}(b)-\nabla\bar{\eta}(b)\cdot(a-b),
 \end{align}
 for all $a,b$ in the state space for the equation \eqref{1-D_Lagrangian}, and similarly, $G(\cdot|\cdot)$ is the relative flux, defined as 
 \begin{align}
     G(a|b)\coloneqq G(a)-G(b)-DG(b)(a-b),
 \end{align}
 $U$ is any  vector of conserved quantities, $U=(\frac{1}{\rho},v_1,v_2)$, which is a bounded, weak solution to \eqref{1-D_Lagrangian} which also verifies the entropy inequality \eqref{entropy_Lagrangian}, and $\bar{U}$ is any classical, Lipschitz solution to \eqref{1-D_Lagrangian}, written as a vector of conserved quantities $\bar{U}=(\frac{1}{\bar{\rho}},\bar{v}_1,\bar{v}_2)$.
 
Furthermore, $\nabla^2 \bar{\eta}$ denotes the Hessian of $\bar{\eta}$ (in Lagrangian variables),
 $s>0$ is a fixed constant which depends on $\bar{\eta},\bar{q}$ and $G$, and $r, t>0$ are any fixed, positive numbers.

 Then, \eqref{key_dafermos} will hold for all points $\sigma$ of $L^\infty$ weak* continuity of $\tau\mapsto\bar{\eta}(U(\cdot, \tau))$ in $(0, t)$.
 
 For more details, see \cite[Equation (5.2.14)]{dafermos_big_book}.

\uline{Step 3}

In the context of \eqref{key_dafermos}, let $U$ be any bounded, weak solution to \eqref{1-D_Lagrangian}, \eqref{entropy_Lagrangian}, depending only $y$ and $t$, with initial data in the form \eqref{init_data} verifying $\rho_+=\rho_-$, $v_{+2}=v_{-2}$.

Remark that given any Lipschitz-continuous function $\bar{v}_1\colon\mathbb{R}\to\mathbb{R}$, a constant $\bar{v}_2\in\mathbb{R}$, and a positive constant $\bar{\rho}>0$, the triple of conserved quantities $(\frac{1}{\bar{\rho}},\bar{v}_1,\bar{v}_2)$ will be a classical solution to \eqref{1-D_Lagrangian}.

Consider then a sequence of smooth functions $\bar{v}_n\colon\mathbb{R}\to\mathbb{R}$ ($n=1,2,\ldots$) such that 
\begin{align}\label{limit_holds}
    \lim_{n\to\infty}\norm{\bar{v}_n(\cdot)-(v_{-1}\mathbbm{1}_{(-\infty,0)}+v_{+1}\mathbbm{1}_{(0,\infty)})}_{L^2(\mathbb{R})}=0.
\end{align}
Such a sequence exists because smooth functions are dense in $L^2$. Define then $\bar{U}_n\coloneqq (\rho_\pm,\bar{v}_n,v_{\pm 2})$, where $\rho_\pm\coloneqq \rho_+=\rho_-$ and similarly for $v_{\pm 2}$.

Consider then \eqref{key_dafermos}, with $U$ as described above, and $\bar{U}_n$ playing the role of $\bar{U}$:

\begin{multline}\label{key_dafermos1}
\int_{|y| < r + s(t - \sigma)} \bar{\eta}(U(y, \sigma)| \bar{U}_n(y, \sigma)) \, dy \leq \int_{|y| < r + st} \bar{\eta}(U(y,0)| \bar{U}_n(y,0)) \, dy
\\
- \int_0^\sigma \int_{|y| < r + s(t - \tau)}  \partial_{y} \bar{U}_n \nabla^2 \bar{\eta}(\bar{U}_n)G(U| \bar{U}_n) \, dy d\tau,
\end{multline}

Notice that in \eqref{key_dafermos1}, the term 
\begin{align}
    \partial_{y} \bar{U}_n \nabla^2 \bar{\eta}(\bar{U}_n)G(U| \bar{U}_n) 
\end{align}
is identically zero, because the Hessian of $\bar{\eta}$ is diagonal, and only the middle component of $\partial_{y} \bar{U}_n$ is nonzero, while the middle component of $G(U| \bar{U}_n)$ is always identically zero.

Moreover, notice that the left-hand size of \eqref{key_dafermos} (and thus also \eqref{key_dafermos1}) is lower-semicontinuous in $\sigma$ due to the convexity of $\bar{\eta}$. Thus, we can take the $\liminf$ as $\sigma\to t^-$ in \eqref{key_dafermos1}. This yields
\begin{align}\label{key_dafermos2}
\int_{|y| < r } \bar{\eta}(U(y,t)| \bar{U}_n(y,t)) \, dy \leq \int_{|y| < r + st} \bar{\eta}(U(y,0)| \bar{U}_n(y,0)) \, dy.
\end{align}

Remark now that for $\lambda_1,\lambda_2,\lambda_3$ and $\lambda_2'$ in state space,
\begin{align}\label{l2_control}
    \bar{\eta}\big((\lambda_1,\lambda_2,\lambda_3)|(\lambda_1,\lambda_2',\lambda_3)\big) = \frac{1}{2}(\lambda_2-\lambda_2')^2,
\end{align}

and similarly, due to the strict convexity of $\bar{\eta}$, we have that $\bar\eta(a|b)=0 \iff a=b$. At this point, we take the limit as $n\to\infty$. Then, from \eqref{limit_holds}, \eqref{key_dafermos2}, and \eqref{l2_control}, we conclude the Theorem.

\appendix

\section{A quantitative inverse function theorem}\label{sec:IFT}

To perturb the approximate solutions to an algebraic system (in our case, given by MATLAB in \Cref{MATLAB_prop}) into exact solutions, we use a quantified Inverse Function Theorem. 

We follow the proof of the very similar quantified inverse function theorem given in \cite[p.~595]{MR0819558}. However, in \cite{MR0819558}, the result is less precise and depends only on the determinant of the Jacobian. Here, we give a more precise  estimate involving the singular values of the Jacobian.

\begin{proposition}[Quantitative Inverse Function Theorem with Singular Values]\label{IFT_prop}

\hfill

Let $F\colon \{ x \in \mathbb{R}^n : |x| \leq 1 \}\to \mathbb{R}^n$ be any $C^2$ function. 

Assume that there exists $r>0$ such that $\sigma_{\text{min}}(DF(0))\geq r$, where $DF$ is the Jacobian of $F$ and given a constant matrix $B$, $\sigma_{\text{min}}(B)$ denotes the smallest singular value of the matrix $B$.

Moreover, assume that there exists $A>0$ such that the second partial derivatives of $F$ are bounded by $A$, i.e. $\abs{\frac{\partial}{\partial x_i x_j} F_k}<A$ for all $i,j,k$ and where $F_k$ is the $k$th component of $F$. 

Then we can conclude that there exists constants $D_1,D_2>0$ (which depend only on $A$ and $n$) such that

\begin{enumerate}
    \item[(i)] $F$ is one-to-one on $\{ |x| \leq D_1 r \}$,
    \item[(ii)] $F( \{ |x| \leq D_1 r \} )$ contains $\{ | y - F(0) | \leq D_2 r^2 \}$.
\end{enumerate}

In particular, we can choose
\begin{align}
D_1 &\coloneqq \frac{1}{4n^2 A}\label{C1_def},\\
D_2 &\coloneqq \frac{1}{8n^2 A}.
\end{align}
\end{proposition}

\begin{figure}[tb]
\centering
      \includegraphics[width=.8\textwidth]{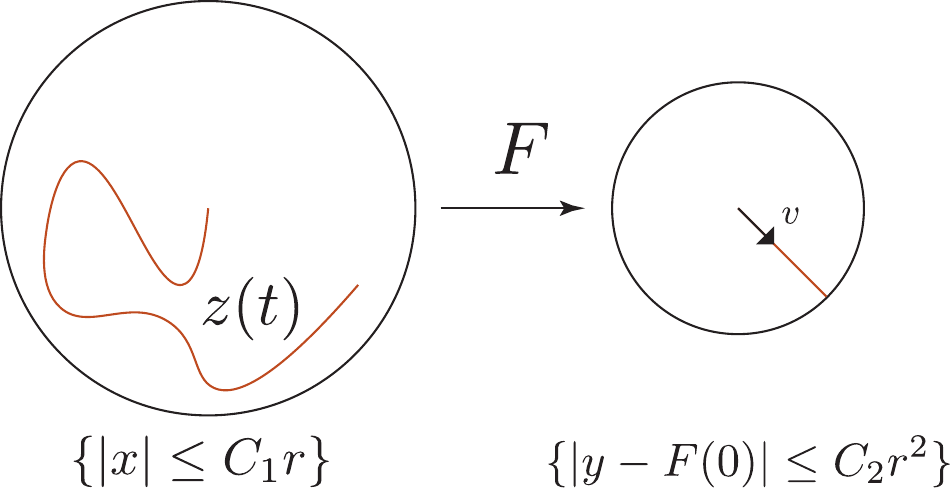}
  \caption{A schematic of the proof of \Cref{IFT_prop}.}\label{fig:IFT}
\end{figure}

\begin{proof}

We recall the following Proposition, which says that singular values of a matrix are 1-Lipschitz in the appropriate norms:

\begin{proposition}[\protect{\cite[Corollary 8.6.2]{MR3024913}}]\label{svd_prop}
    
If $B$ and $B + E$ are in $\mathbb{R}^{m \times n}$ with $m \geq n$, then for $k\in\{1,\ldots,n\}$
\[
| \sigma_k(B + E) - \sigma_k(B) | \leq \sigma_1(E) = \| E \|_2,
\]
where $\sigma_i(B)$ denotes the $i$th largest singular value of the constant matrix $B$.
\end{proposition}

Recall also the matrix norms: the Frobenius norm of a matrix $B\in\mathbb{R}^{m\times n}$
\[
\| B \|_F = \sqrt{\sum_{i=1}^m \sum_{j=1}^n |b_{ij}|^2} 
\]
and the $2$-norm
\[
\| A \|_2 = \sup_{x \neq 0} \frac{\abs{B x}}{\abs{x}},
\]
where $b_{ij}$ denotes the components of $B$.

Note that the matrix $2$-norm is defined in terms of the Euclidean norms of the vectors $Bx$ and $x$.
 
 Thus, from the definition of $D_1$ (see \eqref{C1_def}) and \Cref{svd_prop}, we have that 
 \begin{align}
     \sigma_{\text{min}}(DF(x))\geq \frac{r}{2}
 \end{align}
 for all $x\in\{ x \in \mathbb{R}^n : |x| \leq D_1 r \}$. Remark also that $\| B \|_2\leq\| B \|_F$ for all matrices $B$.

Thus, by the definition of singular values (see \cite[p.~487]{MR3024913}), $|DF(x) \cdot v| \geq \frac{r}{2}|v|$ for all $x$ in this ball and all vectors $v\in\mathbb{R}^n$. For small $x_1$ and $x_2$ given, consider $z(t) \coloneqq F(x_1 + t(x_2 - x_1))$. Then $z'(t) = DF(x_1 + t(x_2 - x_1))(x_2 - x_1)$ and we conclude that  for all $t\in[0,1]$,
\[
\langle z'(t), z'(0) \rangle \geq |z'(0)|^2 - n^2 A |x_2 - x_1|^2 |z'(0)|
\]
\[
\geq (\frac{r}{2} -n^2 A |x_1 - x_2|)|x_1 - x_2||z'(0)|
\]
\[
\geq (\frac{1}{2} - n^2 A D_1) r |x_1 - x_2| \|z'(0)\| > 0
\]
if both $x_i$ live in the set $\{ |x| \leq D_1 r \}$. Remark also our choice of $D_1$ (see \eqref{C1_def}). From this, we receive $[z(1) - z(0), z'(0)] > 0$, and in turn $F(x_1) \neq F(x_2)$. Remark that again we have used that $\| B \|_2\leq\| B \|_F$ for all matrices $B$.

To prove (ii), assume without loss of generality that $F(0) = 0$. Let an arbitrary unit vector $v \in \mathbb{R}^n$ be given. Define the vector field $X$ on $\{ |x| \leq D_1 r \}$ by $(DF(x))X(x) = v$ for all $x$. Then $|X|(x) \leq \frac{2}{r}$ for all $x$ in $\{ |x| \leq D_1 r \}$. Consider the integral curve of $X$ given by the ordinary differential equation $\dot{z}(t) = X(z(t))$ with $z(0)=0$. Notice that $z(t)$ is well-defined for all $0 \leq t \leq \frac{D_1}{2}  r^2$ because of the upper bound on $|X|$ and the existence of solutions to ordinary differential equations. Furthermore, since $\frac{dF(z(t))}{dt} = v$, $F(z(t)) = tv$ for all $t$. Thus (ii) holds for all $D_2 \leq \frac{D_1}{2}$. 

See \Cref{fig:IFT}.

\end{proof}

\section{Explicit values for \Cref{MATLAB_prop}}\label{sec:numerics}

\subsection{Exact values}

Here we give exact values which verify the conclusions of \Cref{MATLAB_prop}. These values are also stored in the comments of our code (available on the GitHub\footnote{See here: \url{https://github.com/sammykrupa/NonUniqueness2DIsentropicEuler}}).

For the purpose of seeing the order of magnitude of the values, and the relations between the various values, we provide decimal approximations of these constants in \Cref{sec:approx}, below.

\begin{align*}
    \hat\alpha_1 &= -\frac{8177336068870495}{140737488355328}\\
    \hat\alpha_2 &=-\frac{4833381446756075}{562949953421312}\\
    \hat\alpha_3 &= \frac{3121572020159473}{562949953421312}\\
    \hat\beta_1 &= -\frac{2536561643647751}{140737488355328}\\
    \hat\beta_2 &=-\frac{1114286601116939}{70368744177664}\\
    \hat\beta_3 &=-\frac{6219197795695073}{562949953421312}\\
    \hat\gamma_1 &=\frac{841617150350781}{549755813888}\\
    \hat\gamma_2 &=-\frac{2850833975067331}{17592186044416}\\
    \hat\gamma_3 &=-\frac{8954832877447991}{140737488355328}\\
    \hat\delta_1 &=\frac{28872176135415855785280523654056524019908546591}{27627619078169805047324756605549438692229120}\\
    \hat\delta_2 &=-\frac{867454945412067709200232997995952542374584537720982074207241594308599438377}{32748846874784971211058574285222379723549486466626634273206947431338475520}\\
    \hat\delta_3 &=-\frac{2871256077954219}{35184372088832}\\
    \hat\rho_- &=\hat\rho_+ =\frac{2708112612978501}{281474976710656}\\
    \hat\rho_1 &=\frac{6811063536043807}{562949953421312}\\
    \hat\rho_2 &=\frac{2057060350258899}{562949953421312}\\
    \hat\rho_3 &=\frac{3062207031116133}{281474976710656}\\
    \hat\epsilon_- &=\hat\epsilon_+=-\frac{5041529442624971}{2199023255552}\\
    \hat\epsilon_1 &=-\frac{5015532875605977}{2199023255552}\\
    \hat\epsilon_2 &=-\frac{5073206593829053}{2199023255552}\\
    \hat\epsilon_3 &=-\frac{2515400677054201}{1099511627776}\\
    \hat\epsilon_-' &=\hat\epsilon_+' =\frac{1676289422169645}{562949953421312}\\
    \hat\epsilon_1' &=\frac{60006068216738166351756926651195471316209797920723479281182349}{9106752347169134708406278810788569756749285046122185710632960}\\
    \hat\epsilon_2' &=\frac{1831278218949756891087541424381121781924211556364089293813566516592411003}{1359848686641341079894728790850171965963388817622134940140753492461486080}\\
    \hat\epsilon_3' &=\frac{5400383921383283}{1125899906842624}\\
    \hat\nu_- &=-\frac{6486283176597739958874052307549}{196306040423407104692364247040}\\
    \hat\nu_1 &=-\frac{1153852086001065889673487658885}{60824224363690518566334889984}\\
    \hat\nu_2 &=-\frac{4856156003780791}{562949953421312}\\
    \hat\nu_+ &=\frac{7162856387903725}{562949953421312}\\
    \hat v_{-1} &=-\frac{4098844157247653}{70368744177664}\\
    \hat v_{+1} &=\frac{3603433899522037}{562949953421312}\\
    \hat v_{-2} &= \hat v_{+2} =-\frac{996118042660627}{70368744177664}\\
    \hat C_1 &= \frac{510415269881361}{137438953472}\\
    \hat C_2 &=\frac{1515879700153707}{2199023255552}\\
    \hat C_3 &=\frac{1855257252703141}{8796093022208}
    \end{align*}

\subsection{Approximate values}\label{sec:approx}

\begin{align*}
    \hat\alpha_1 &\approx -58.1,\quad \hat\alpha_2 \approx-8.59,\quad \hat\alpha_3 \approx 5.55\\
    \hat\beta_1 &\approx -18.0,\quad \hat\beta_2 \approx-15.8,\quad \hat\beta_3 \approx-11.0\\
    \hat\gamma_1 &\approx1530.0,\quad \hat\gamma_2 \approx-162.0,\quad \hat\gamma_3 \approx-63.6\\
    \hat\delta_1 &\approx1050.0,\quad \hat\delta_2 \approx-26.5,\quad \hat\delta_3 \approx-81.6\\
    \hat\rho_-&=\hat\rho_+ \approx9.62,\quad \hat\rho_1 \approx12.1,\quad \hat\rho_2 \approx3.65,\quad \hat\rho_3 \approx10.9\\
    \hat\epsilon_-&= \hat\epsilon_+\approx-2290.0,\quad \hat\epsilon_1 \approx-2280.0,\quad \hat\epsilon_2 \approx-2310.0,\quad \hat\epsilon_3 \approx-2290.0\\
    \hat\epsilon_-' &=\hat\epsilon_+' \approx2.98,\quad \hat\epsilon_1' \approx6.59,\quad \hat\epsilon_2' \approx1.35,\quad \hat\epsilon_3' \approx4.8\\
    \hat\nu_- &\approx-33.0,\quad \hat\nu_1 \approx-19.0,\quad \hat\nu_2 \approx-8.63,\quad \hat\nu_+ \approx12.7\\
    \hat v_{-1} &\approx-58.2,\quad \hat v_{+1} \approx6.4,\quad \hat v_{-2}=\hat v_{+2} \approx-14.2\\
    \hat C_1 &\approx 3710.0,\quad \hat C_2 \approx689.0,\quad \hat C_3 \approx211.0
    \end{align*}

As a final remark we note that the values of the internal energy $\epsilon(\rho_i)=\hat\epsilon_i$ and $\epsilon(\rho_{\pm})=\hat\epsilon_{\pm}$ are negative. However, using the continuity equation it is easy to see that one can add an arbitrary constant to the internal energy $\epsilon(\rho)$ without changing the property of being a weak solution/subsolution, changing the pressure law \eqref{pressure_relation} or changing the inequalities \eqref{convex_ineq} ensuring strict hyperbolicity.

\bibliographystyle{plain}
\bibliography{references}
\end{document}